\documentclass[twoside,11pt]{article}

\usepackage{amsfonts, amsmath, amssymb, amsthm, constants}
\usepackage{graphicx, subfigure}
\usepackage[top=1in,bottom=1in,left=1in,right=1in]{geometry}
\usepackage[sort&compress]{natbib} 

\usepackage{color}
\definecolor{darkred}{RGB}{150,0,0}
\definecolor{darkgreen}{RGB}{0,150,0}
\definecolor{darkblue}{RGB}{0,0,200}

\usepackage{hyperref}
\hypersetup{colorlinks=true, linkcolor=darkred, citecolor=darkgreen, urlcolor=darkblue}
\usepackage{url}

\newtheorem{thm}{Theorem}

\newtheorem{lem}{Lemma}

\def\beq{\begin{equation}}
\def\eeq{\end{equation}}
\def\beqn{\begin{eqnarray*}}
\def\eeqn{\end{eqnarray*}}

\newcommand{\thmref}[1]{Theorem~\ref{thm:#1}}

\newcommand{\lemref}[1]{Lemma~\ref{lem:#1}}
\newcommand{\secref}[1]{Section~\ref{sec:#1}}

\DeclareMathOperator*{\argmin}{arg\, min}

\def\cB{\mathcal{B}}

\def\cK{\mathcal{K}}

\def\cR{\mathcal{R}}

\def\bbN{\mathbb{N}}

\def\bbR{\mathbb{R}}

\def\bbT{\mathbb{T}}

\def\bbZ{\mathbb{Z}}

\newcommand{\expect}[1]{\mathbb{E}\left(#1\right)}
\newcommand{\pr}[1]{\mathbb{P}\left(#1\right)}

\def\eps{\varepsilon}

\def\iid{\stackrel{\rm i.i.d.}{\sim}}
\def\Bin{\text{Bin}}

\def\cst{C}

\def\cat{{\rm Cat}}
\def\gen{{\rm gen}}

\def\rA{\mathring{A}}
\def\rB{\mathring{B}}
\def\uA{\underline{A}}
\def\uB{\underline{B}}

\begin{document}

\title{Finite Size Percolation in Regular Trees}
\author{ 
Ery Arias-Castro \\[.1in]
University of California, San Diego
\footnote{
Department of Mathematics, University of California, San Diego
\{\href{http://math.ucsd.edu/~eariasca}{http://math.ucsd.edu/{\tt \~{}}eariasca}\}
}
}
\date{}
\maketitle

\begin{abstract}
In the context of percolation in a regular tree, we study the size of the largest cluster and the length of the longest run starting within the first $d$ generations.  As $d$ tends to infinity, we prove almost sure and weak convergence results.

\medskip

\noindent {\bf MSC 2010:} Primary 60K35.
 
\noindent {\bf Keywords and Phrases:}  Percolation on trees; largest open cluster; longest success run; Galton-Watson processes; Chen-Stein method for Poisson approximation.
\end{abstract}

\section{Introduction}
Fix a positive integer $r$ and let $\bbT$ be the infinite $r$-ary tree, rooted at $\rho_0$.  
We consider a Bernoulli percolation on $\bbT$. 
Formally, to each node $v \in \bbT$, we associate a random variable $X_v$, where the variables $\{X_v:v \in \bbT\}$ are i.i.d.~Bernoulli with $\pr{X_v = 1} = 1 -\pr{X_v = 0} = p \in (0,1)$.  
For a subset $A \subset \bbT$, let $X_A = \prod_{v \in A} X_v$.  We say that $A$ is {\it open} if $X_A = 1$.

\subsection{The size of the largest cluster.}
We use the term {\it cluster} to denote a connected component (i.e.~subtree) of $\bbT$ when undirected.  Let $\cK$ denote the set of clusters in $\bbT$.
For a node $v \in \bbT$, let $\gen(v)$ be its generation, i.e.~the number of nodes in the shortest path from the root $\rho_0$ to $v$, not counting $\rho_0$.  Note that $\gen(\rho_0) = 0$.
Let $\bbT_d$ be the set of nodes with generation not exceeding $d$, namely $\bbT_d = \{v \in \bbT: \gen(v) \leq d\}$.
For a cluster $A \in \cK$, we let $|A|$ denote its size (i.e.~number of nodes) and $\rho(A)$ its root, namely $\rho(A) = \argmin\{\gen(v): v \in A\}$.
For $d \in \bbN$, define $K_d$ to be the size of the largest open cluster with root of generation not exceeding $d$:
$$K_d = \max\{|A|: A \in \cK,\, \rho(A) \in \bbT_d,\, X_A = 1\}.$$
In particular, $K_0$ is the size of the largest open cluster containing the root $\rho_0$.

In this paper we study the limit behavior of $K_d$, as $d \to \infty$.
In the context of the one-dimensional lattice $\bbZ$, the corresponding results are often referred to as the Erd\"os-R\'enyi Law~\cite{MR0272026} and, in that context, our approach follows that of Arratia, Goldstein and Gordon~\cite{MR972770}.   
In higher dimensions, the problem is much more intricate and many questions remain without answer.  For a sample of sophisticated results, see e.g.~\cite{MR1868996,MR1372330,MR1880230}.  The book by Grimmett~\cite{MR1707339} is a standard reference on percolation.  For references more specific to trees, we refer the reader to a survey paper by Pemantle~\cite{MR1368099} and the book of Lyons and Peres~\cite{tree-book}.  Though the literature on percolation is vast, most of it focuses on the existence of an infinite cluster and its characteristics when it exists.
On the applications side, Patil and Taillie~\cite{MR2109372} identify regions of interest in a network by thresholding the response from each site in the network and computing connected components, which amounts to extracting the open clusters.  One imagines that the largest cluster might receive the most attention.  In particular, they mention monitoring water quality in a network of freshwater streams, where each stream may be modeled as a tree, though an irregular one. 
%

It is well-known that, in the supercritical setting where $p > 1/r$, the cluster at the origin has positive probability of being infinite, and, in fact, $\pr{K_d = \infty}$ tends to 1 as $d$ increases.
We restrict our attention to the subcritical and critical cases, i.e~$p < 1/r$ and $p = 1/r$ respectively, where the cluster at the origin is finite with probability one.
We start with the critical case, where we show that $K_d$ behaves like the maximum of $r^d$ independent random variables with distribution the total progeny of a Galton-Watson process with offspring distribution $\Bin(r, 1/r)$.  Let $\log_r$ denote the logarithm in base $r$.
\begin{thm} \label{thm:cluster-cri}
Assume $p=1/r$.  Then with probability one,
$$\frac{\log_r K_d}{d} \rightarrow 2, \quad d \to \infty.$$
Moreover,
\[
\pr{\log_r K_d \leq 2d +x} \to \to \exp(- \Cl{cluster-cri}\, r^{-x/2}), \quad d \to \infty, \quad \Cr{cluster-cri} := \frac{2\, n^{-1/2}}{\sqrt{2 \pi r (r-1)}}
\]
\end{thm}
%

For the subcritical case, we obtain similar results without transformation.  Here, a Poisson approximation applies showing that $K_d$ behaves like the maximum of $|\bbT_d| = r^{d+1}/(r-1)$ independent random variables with distribution the total progeny of a Galton-Watson process with offspring distribution $\Bin(r, p)$.
Define 
\begin{equation}
\label{kappa}
\kappa = p (1-p)^{r-1}\ \frac{r^r}{(r-1)^{r-1}}.
\end{equation}
Note that $\kappa < 1$ for all $p < 1/r$.  Let $[x]$ denote the entire part of $x \in \bbR$.
\begin{thm}
\label{thm:cluster-sub}
Assume $p<1/r$.  Then with probability one,
$$\frac{K_d}{d} \rightarrow \frac{1}{\log_r(1/ \kappa)}, \quad d \to \infty.$$
Moreover, the sequence of random variables $(K_d -\mu_d: d \geq 0)$ is tight, where $\mu_d := \frac{d -\frac{3}{2} \log_r d}{\log_r(1/ \kappa)}$.  In addition, a subsequence $K_d -\mu_d$ converges weakly if, and only if, $a := \lim_{d \to \infty} (\mu_d -[\mu_d])$ exists, in which case the weak limit is $[Z+a] -a$, where
\[
\pr{Z \leq z} = \exp\left(-\Cl{cluster}\, \kappa^{z}\right),
\]
for an explicit constant $\Cr{cluster} > 0$ depending only on $(p,r)$
\end{thm}
The behavior of the size of the largest open cluster in the subcritical regime is therefore similar in the context of the regular tree and in the context of the one-dimensional lattice, the latter corresponding to the length of the longest perfect head run in a sequence of coin tosses~\cite[Ex.~3]{MR972770}.

\subsection{The length of the longest run.}
We use the term {\it run} for a path in $\bbT$ when directed away from the root $\rho_0$.
Note that runs are special clusters. 
Let $\cR$ denote the set of runs and define $R_d$ to be the length of the longest open run with root of generation not exceeding $d$:
$$R_d = \max\{|A|: A \in \cR,\, \rho(A) \in \bbT_d,\, X_A = 1\}.$$

Of course, runs and clusters coincide in the one-dimensional lattice $\bbZ$.  For a general reference on runs in dimension one, see~\cite{MR1882476}.  Using the Chen-Stein method, Chen and Huo~\cite{MR2268049} proved results on the longest left-right run in a thin two-dimensional lattice of the form $([0,d]\times[0,a]) \cap \bbZ^2$, with the width $a$ remaining constant.  We also mention the work of Arias-Castro, Donoho and Huo~\cite{MR2275244} who used a statistic based on the longest run in a particular, non-planar graph to detect filaments in point-clouds.

The results we obtain for runs are parallel to those we obtain for clusters.
In the critical case, we show that $R_d$ behaves like the maximum of $r^d$ independent random variables with distribution the height of a Galton-Watson process with offspring distribution $\Bin(r, 1/r)$. 

\begin{thm} \label{thm:run-cri}
Assume $p=1/r$.  Then with probability one,
$$\frac{\log_r R_d}{d} \rightarrow 1, \quad d \to \infty.$$
Moreover, for any $x \in \bbR$,
\[
\pr{\log_r R_d \leq d + x}  \to \exp(- \Cl{run-cri}\, r^{-x}), \quad d \to \infty, \quad \Cr{run-cri} := \frac{2 rp}{r-1}.
\]  
\end{thm}

In the subcritical case, we show that $R_d$ behaves like the maximum of $|\bbT_d|$ independent random variables with distribution the height of a Galton-Watson process with offspring distribution $\Bin(r, p)$.  Again, a Poisson approximation applies.  The constant that appears in the exponent is only defined implicitly.
\begin{thm}
\label{thm:run-sub}
Assume $p<1/r$.  Then with probability one,
$$\frac{R_d}{d} \rightarrow \frac{1}{\log_r (1/p) -1}, \quad {\rm as} \ d
\rightarrow \infty.$$
Moreover, the sequence of random variables $(K_d -\nu_d: d \geq 0)$ is tight, where $\nu_d := \frac{d}{\log_r (1/p) -1}$.  In addition, a subsequence $K_d -\nu_d$ converges weakly if, and only if, $a := \lim_{d \to \infty} (\nu_d -[\nu_d])$ exists, in which case the weak limit is $[Z+a] -a$, where
\[
\pr{Z \leq z} = \exp\left(-\Cl{run}\, (rp)^{z}\right),
\]
for an explicit constant $\Cr{run} > 0$ depending only on $(p,r)$
\end{thm}

\subsection{Contents.}
The rest of the paper is devoted to proving our results.  In \secref{proof-cluster} we prove \thmref{cluster-cri} and \thmref{cluster-sub}.  In \secref{proof-run} we prove \thmref{run-cri} and \thmref{run-sub}.

\subsection{Additional Notation.}
Let $\partial \bbT_d = \{v \in \bbT: \gen(v) = d\}$.
For a cluster $A$, let $\uA$ denote the set of nodes not in $A$ whose parents belong to $A$, and if $\rho(A) \neq \rho_0$, let $\rA$ denote the parent of $\rho(A)$.
Also, define $(1 -X)_A = \prod_{v \in A} (1 - X_v)$.
For two sequences of real numbers $(a_n)$ and $(b_n)$, we use the notation $a_n \sim b_n$ to indicate that $a_n/b_n \to 1$ and $a_n \asymp b_n$ to indicate that the ratio $a_n/b_n$ is bounded away from zero and infinity, both understood as $n \to \infty$.
Throughout the paper $C$ denotes a finite, positive constant depending only on $r$ and $p$, whose value may change with each appearance.


\section{The size of the largest open cluster}
\label{sec:proof-cluster}

In this section, we prove \thmref{cluster-cri} and \thmref{cluster-sub}.
We start with some notation.
For a vertex $v \in \bbT$, let $K(v)$ be the size of the largest open cluster with root $v$, 
$$K(v) = \max\{|A|: A \in \cK,\, \rho(A) = v,\, X_A = 1\}.$$
In particular,
$$K_d = \max\{K(v): v \in \bbT_d\}.$$
The distribution of $K(v)$ does not depend on $v \in \bbT$, and, in fact, given $X_v = 1$, coincides with that of the total progeny of a Galton-Watson tree starting with one individual and with offspring distribution Bin$(r,p)$.
Define 
\[
\psi_n = \pr{K(v) = n}, \quad \Psi_n = \pr{K(v) > n}.
\]
Applying a well-known identity by Dwass~\cite{MR0253433} (called the Otter-Dwass formula in~\cite{tree-book}), we get
\beqn
\psi_n 
&=& \frac{p}{n} \pr{\xi_1 + \cdots + \xi_n = n-1}, \text{ where } \xi_1, \dots, \xi_n \iid \Bin(r, p) \\
&=& \frac{p}{n} \pr{\Bin(n r, p) = n-1} \\
&=& \cat_{n}\, p^{n} (1-p)^{n (r-1) + 1},
\eeqn
where
\[
\cat_{n} := \frac{1}{n} {n r \choose n-1} = \frac{1}{(r-1) n + 1} {r n \choose n} 
\]
is the {\it $n$th generalized Catalan number}~\cite{hil91}, which among other interpretations, is the number of subtrees of $\bbT$ of size $n$ rooted at the origin, i.e.~
\[
\cat_{n} = |\{A \in \cK: \rho(A) = \rho_0,\, |A| = n\}|.
\]
We could have obtained the expression for $\psi_n$ using this definition of $\cat_n$. Indeed, for $n  > 0$, $K(v) = n$ if, and only if, there is a (unique) subtree $A$ with $|A| = n$ , $\rho(A) = v$ and $X_A (1 -X)_{\uA} = 1$, so that $A$ cannot be extended and still be an open cluster.  We then use the fact that a subtree of size $n$ has exactly $(r-1)n+1$ children.
With the use of Stirling's formula, we arrive at the following conclusions; see also~\cite{MR0386042, ney}. 

\begin{lem}
\label{lem:Kv}
In the critical case $p  =1/r$,
\[
\Psi_n \sim \frac{\Cr{cluster-cri}}{\sqrt{n}}.
\]
In the subcritical case $p < 1/r$, 
\[
\Psi_n \sim \Cl{cluster-aux}\, \frac{\kappa^{n+1}}{n^{3/2}}, \quad \Cr{cluster-aux} := \frac{1}{\sqrt{2 \pi}(1-\kappa)} \frac{(1-p) r^{1/2}}{(r-1)^{3/2}}.
\]
\end{lem}

\subsection{Proof of \thmref{cluster-cri}}
\label{sec:proof-cluster-cri}

Define
\[
K^\partial_d := \max\{K(v): v \in \partial \bbT_d\}.
\]
We first prove that the conclusions of \thmref{cluster-cri} hold for $K^\partial_d$.  
For $x \in \bbR$, let $n_d(x) = [r^{2 d + x}].$
As $K^\partial_d$ only involves independent random variables, we have
\begin{equation} \nonumber 
\pr{\log_r K^\partial_d \leq 2d +x} = \pr{K^\partial_d \leq n_d(x)} = (1 -\Psi_{n_d(x)})^{r^d} = \exp(- \Cr{cluster-cri}\, r^{-x/2} + O(r^{-d -x})).
\end{equation}
Letting $d \to \infty$, we obtain the weak convergence, and by choosing $x = \eps d$, with $\eps > -2$ fixed, and applying the Borel-Cantelli Lemma, we obtain the almost sure convergence.

It therefore suffices to show that $K_d = (1 + o_P(1)) K^\partial_d$.  Clearly, $K_d \geq K^\partial_d$, so we focus on the upper bound.    
Define 
\[
B_d = \{v \in \partial \bbT_d: K(v) > r^d/d\}, \quad B^2_d = \{v \in \partial \bbT_d: K(v) > r^{2d}/d\}.
\]
For any open cluster $A$ with $\rho(A) \in \bbT_d$, we have
\[
|A| = |A \cap \bbT_{d-1}| + \sum_{v \in A \cap \partial \bbT_d} K(v) \leq r^d + r^{2d}/d + \sum_{v \in A \cap B_d} K(v).
\]
We turn to bounding the sum.  We first show that, with probability tending to one, there is no open cluster $A$ containing three or more nodes in $B_d$.  Indeed, take $v_1, v_2, v_3 \in \partial \bbT_d$ distinct.  Let $w$ denote their most recent common ancestor and let $k = d -\gen(w)$.  Either the paths $v_j \to \rho_0$ meet at $w$ for the first time or two of the paths meet at a node $u$ with $\gen(u) > \gen(w)$, in which case we let $\ell = d -\gen(u)$.  Now, the nodes $v_1, v_2, v_3$ belong to the same open cluster if, and only if, the smallest subtree containing $w$ and $v_1, v_2, v_3$ is open, and this subtree is of size $\ell + 2 k + 1$, and therefore, the probability that they belong to the same open cluster is $p^{\ell + 2k + 1}$.  In addition, the number of such triplets is bounded by 
\[
\left(r^{d-k} {r^k \choose 3}\right) \cdot \left(r^k r^{k -\ell} {r^\ell \choose 2}\right) {r^k \choose 3}^{-1} \asymp r^{d +k +\ell}. 
\] 
The first factor comes from the fact that the three nodes are leaves of a subtree with root at generation $d-k$.  Given that, the second factor comes from the fact that two of them belong to a subtree of that subtree with root at (relative) generation $k -\ell$.  
Hence, remembering that $p = 1/r$ and using \lemref{Kv}, we have
\beqn
\pr{\exists A \in \cK: X_A = 1,\, |A \cap B_d| \geq 3}
&\leq& C\, \pr{K(v) > r^d/d}^3 \cdot \sum_{k=0}^d \sum_{\ell=0}^k r^{d +k +\ell} p^{\ell + 2k + 1} \\
&\leq& C\, (r^d/d)^{-3/2} r^{d} \asymp d^{3/2} r^{-d/2}.
\eeqn
By the same token, with probability tending to one (in fact of order at most $d/r^d$), there is no open cluster $A$ containing two or more nodes in $B^2_d$.
Now, when $|A \cap B_d| \leq 2$ and $|A \cap B^2_d| \leq 1$, we have
\[
\sum_{v \in A \cap B_d} K(v) \leq \max_{v \in A \cap B_d} K(v) + r^{2d}/d \leq K^\partial_d + r^{2d}/d.
\]
In the end, with probability tending to one, 
\[
|A| \leq r^d + 2\, r^{2d}/d + K^\partial_d,
\]
for any open cluster $A$ with $\rho(A) \in \bbT_d$.  Hence,
\[
K_d \leq K^\partial_d + O_P(r^{2d}/d),
\]
and we conclude by the fact that $K^\partial_d$ is of order exceeding $r^{2d}/d$ with probability tending to one.

\subsection{Proof of \thmref{cluster-sub}}
\label{sec:proof-cluster-sub}

The proof of the almost sure convergence may be obtained following the arguments provided in \secref{proof-cluster-cri} or using the bounds we are about to prove below.  We omit details.


The proof of the weak convergence is based on the Chen-Stein method for Poisson approximation as formulated by Arratia, Goldstein and Gordon~\cite{MR972770}.  
Define
$$Y_A = \left\{\begin{array}{ll}
X_A (1 -X)_{\uA}, & \rho(A) = \rho_0,\\[.05in]
X_A (1-X)_{\rA} (1 -X)_{\uA}, & \rho(A) \neq \rho_0;
\end{array}\right.$$
Also, let $\cK_{d,n}$ be the set of clusters of size exceeding $n$ with root in $\bbT_d$, and define
$$W_{d,n} = \sum_{A \in \cK_{d,n}} Y_A.$$
\vspace{-.1in}
By definition,
$$\{K_d \leq n\} = \{Y_A = 0,\, \forall A \in \cK_{d,n}\} = \{W_{d,n} = 0\}.$$
We approximate the law of $W_{d,n}$ by the Poisson distribution with same mean $\lambda_{d,n} = \expect{W_{d,n}}$.  
We start by estimating $\lambda_{d,n}$ using \lemref{Kv}, obtaining
\begin{eqnarray*}
\lambda_{d,n} 
& = & \sum_{A \in \cK_{d,n}} \pr{Y_A = 1} \\
& = & \pr{K(\rho_0) > n} + (1-p) \sum_{v \in \bbT_d, v \neq \rho_0} \pr{K(v) > n} \\
& = & \Psi_n + (1-p) (|\bbT_d| -1) \Psi_n.
\end{eqnarray*}
In particular, as $n, d \to \infty$, 
\[
\lambda_{d,n} \sim \Cr{cluster}\ r^d n^{-3/2} \kappa^{n+1}, \quad \Cr{cluster} := \frac{\Cr{cluster-aux} (1-p) r}{r-1}.
\]
 
For a cluster $A \in \cK_{d,n}$, define its neighborhood $\cB(A)$ as the set of clusters $B \in \cK_{d,n}$ such that 
$$(\rB \cup B \cup \uB) \ \cap \ (\rA \cup A \cup \uA) \neq \emptyset.$$
Define the following sums
\begin{eqnarray*}
F_{d,n} & = & \sum_{A \in \cK_{d,n}} \sum_{B \in \cB(A)} \pr{Y_A = 1} \pr{Y_B = 1},\\
G_{d,n} & = & \sum_{A \in \cK_{d,n}} \sum_{B \in \cB(A), B \neq A} \pr{Y_A = Y_B = 1},\\
H_{d,n} & = & \sum_{A \in \cK_{d,n}} \expect{\left|\expect{Y_A - \expect{Y_A}|Y_B, B \notin \cB(A)}\right|}.
\end{eqnarray*}
Then by the second part of~\cite[Th.~1]{MR972770}, 
$$\left|\pr{W_{d,n} = 0} - \exp(-\lambda_{d,n})\right| \leq F_{d,n} + G_{d,n} + H_{d,n}.$$
For $x \in \bbR$, define $n_d(x) = \left[\mu_d +x\right]$.
When $x$ is fixed and $d \to \infty$, $\lambda_{d,n_d(x)} \asymp 1$, with 
$$\lambda_{d,n_d(x)} \to \Cr{cluster}\, \kappa^{[a+x] -a +1}, \ \text{ when } \mu_d -[\mu_d] \to a, \ x -[x] \neq 1-a,$$
with 
$$\pr{[Z+a] -a \leq x} = \exp(-\Cr{cluster}\, \kappa^{[a+x] -a +1}).$$
Therefore, to conclude it suffices to prove that $F_{d,n}, G_{d,n}, H_{d,n} \to 0$ when $d,n \to \infty$ in such a way that $\lambda_{d,n} \asymp |\bbT_d| \Psi_n \asymp 1$.
First, $H_{d,n} = 0$ by independence of $Y_A$ and $Y_B, B \notin \cB(A)$.  
For $G_{d,n}$, the only pairs $A,B \in \cK_{d,n}$ that contribute to the sum satisfy either $\rB \in \uA$ or $\rA \in \uB$, and in both cases
$$\pr{Y_A = Y_B = 1} = (1-p)^{-1} \pr{Y_A = 1} \pr{Y_B = 1}.$$
Hence, using the fact that there are $\cat_{m}$ subtrees of size $m$ with a given root, each with $(r-1)m + 1$ children, and then \lemref{Kv}, we have
\begin{eqnarray*}
G_{d,n} 
& \leq & 2 (1-p)^{-1} |\bbT_d| \sum_{m > n}\ \cat_{m}\, p^m (1-p)^{(r-1)m+1} ((r-1)m + 1) \cdot \Psi_n \\
& \leq & \cst\, \lambda\, \sum_{m > n} m \psi_{m} = \cst\, \lambda\, \left((n+1) \Psi_n + \sum_{m > n} \Psi_{m}\right)
 \asymp n^{-1/2} \kappa^n \to 0, \quad n \to \infty.
\end{eqnarray*}
For $F_{d,n}$, the only pairs $A,B \in \cK_{d,n}$ that contribute to the sum satisfy either $\rB \in \rA \cup A \cup \uA$ or $\rA \in \rB \cup B \cup \uB$.  The computations are then similar.

\section{The length of the longest open run}
\label{sec:proof-run}
The arguments are parallel to those provided in \secref{proof-cluster}.
For $A \subset \bbT$, define its height as $\tau(A) = \sup\{\gen(v): v \in A\} -\gen(\rho(A))$.
For a vertex $v \in \bbT$, let $R(v)$ be the length of the longest run with root $v$, 
\[
R(v) = 1 + \max\{\tau(A): A \in \cK, \rho(A) = v\}.
\]
In particular, 
\[
R_d = \max\{R(v): v \in \bbT_d\}.
\]
The distribution of $R(v)$ does not depend on $v \in \bbT$, and, in fact, given $X_v = 1$, coincides with that of the height (plus one), i.e.~extinction time, of a Galton-Watson tree with offspring distribution Bin$(r,p)$.
Define 
\[
\phi_h = \pr{R(v) = h}, \quad \Phi_h = \pr{R(v) > h}.
\]
We have the following results on the asymptotic behavior of $\Phi_h$~\cite{ney}.
  
\begin{lem}
\label{lem:Rv}
In the critical case $p  =1/r$, 
\[
\Phi_h \sim \frac{\Cr{run-cri}}{h}.
\]
In the subcritical case $p < 1/r$, there is an implicit constant $\Cl{run-aux} > 0$ such that 
\[
\Phi_h \sim \Cr{run-aux}\, (rp)^h.
\]
\end{lem}

Let $\cat_{n,h}$ denote the number of subtrees rooted at the origin, of size $n$ and height $h$.  See~\cite{MR1249127} for some results on $\cat_{n,h}$.
As in \secref{proof-cluster}, we can argue that 
\[
\phi_h = \sum_{n > h} \cat_{n,h}\, p^n (1-p)^{(r-1)n +1}, \ \text{ implying } \ \Phi_h = \sum_{\ell > h} \sum_{n > \ell} \cat_{n,\ell}\, p^n (1-p)^{(r-1)n +1}.
\]

\subsection{Proof of \thmref{run-cri}}

The proof is based on the following observation
\[
R^\partial_d \leq R_d \leq R^\partial_d + d, \quad R^\partial_d := \max\{R(v): v \in \partial \bbT_d\}, 
\]
where the $d$ term bounds the length of any run in $\bbT_{d-1}$.
As $R^\partial_d$ only involves independent random variables, 
\begin{equation} \label{eq:partial-Rd}
\pr{R^\partial_d \leq h} = (1 -\Phi_{h})^{r^d}. 
\end{equation}
Choosing $h = r^{[d/2]}$ and using \lemref{Rv}, we obtain
\[
\pr{R^\partial_d \leq r^{[d/2]}} \leq \exp(- C\, r^{-d/2}), \ \text{ for some } C > 0,
\]
so that, applying the Borel-Cantelli Lemma, $R^\partial_d \geq r^{d/2}$ eventually, with probability one.  Hence, $\log_r R_d = (1 +o_P(1)) \log_r R^\partial_d$, and it is therefore enough to prove the results for $R^\partial_d$ in place of $R_d$.  
The almost sure convergence is obtained in a similar way by choosing $h = r^{(1+\eps) d}$ with $\eps$ fixed, either positive or negative. 
For the weak convergence, fix $x$ and let $h_d(x) = [r^{d+x}]$.  By \lemref{Rv} and \eqref{eq:partial-Rd}, we have 
\[
\pr{\log_r R^\partial_d \leq d + x}  \to \exp(- \Cr{run-cri}\, r^{-x}), \quad d \to \infty.
\]

\subsection{Proof of \thmref{run-sub}}

We again omit the details of the proof of the almost sure convergence and focus on proving the weak convergence.  Let $\cK_{d,h}$ denote the set of clusters with root in $\bbT_d$ and height exceeding $h$.  We use the notation introduced in \secref{proof-cluster-sub}, with $\cK_{d,h}$ in place of $\cK_{d,n}$.
%
By definition,
$$\{R_d \leq h\} = \{W_{d,h} = 0\}.$$
Using \lemref{Rv}, we obtain
\begin{eqnarray*}
\lambda_{d,h} 
& = & \sum_{A \in \cK_{d,h}} \pr{Y_A = 1}\\
& = & \pr{R(\rho_0) > h} + (1-p) \sum_{v \in \bbT_d} \pr{R(v) > h}\\ 
& = & \Phi_h + (1-p) (|\bbT_d| -1) \Phi_h.
\end{eqnarray*}
In particular, as $h, d \to \infty$, 
\[
\lambda_{d,h} \sim \Cr{run}\, r^d (rp)^{h+1}, \quad \Cr{run} := \frac{\Cr{run-aux}\, (1-p)}{p(r-1)}.
\]

%
For $x \in \bbR$, define $h_d(x) = [\nu_d + x]$.  When $x$ is fixed and $d \to \infty$, we have $\lambda_{d,h_d(x)} \asymp 1$, with   
$$\lambda_{d,h_d(x)} \to \Cr{run} (rp)^{[a+x] -a +1}, \ \text{ when } \nu_d -[\nu_d] \to a, \ x -[x] \neq 1-a.$$
It then suffices to show that $F_{d,h}, G_{d,h}, H_{d,h} \to 0$ when $d,h \to \infty$ in such a way that $\lambda_{d,h} \asymp |\bbT_d| \Phi_h \asymp 1$, and the computations are parallel to those in \secref{proof-cluster-sub}.  
%
We focus on $G_{d,h}$.  Fix $\tilde{p} \in (p, 1/r)$ and let $\tilde{\Phi}_h$ be defined as $\Phi_h$, with $\tilde{p}$ in place of $p$.  For $h$ large enough, we then have
\begin{eqnarray*}
G_{d,h} 
&\leq& 2 \sum_{A \in \cK_{d,h}}\ \sum_{\cB(A), \rB \in \uA} (1-p)^{-1} \pr{Y_A = 1} \pr{Y_B = 1} \\
& \leq & \cst\, |\bbT_d| \sum_{\ell > h} \sum_{n > \ell} \cat_{n,\ell}\ p^n (1-p)^{(r-1)n + 1} ((r-1)n+1) \cdot  \Phi_h \\
& \leq & \cst\, \lambda\, \sum_{\ell > h} \sum_{n > \ell} \cat_{n,\ell}\ \tilde{p}^n (1-\tilde{p})^{(r-1)n + 1} \\
& = & \cst\, \lambda\, \tilde{\Phi}_h  
 \asymp (r \tilde{p})^h \to 0, \quad d \to \infty. 
\end{eqnarray*}
%

\subsection*{Acknowledgements}
The author would like to thank Philippe Flajolet for fruitful conversations and Jason Schweinsberg for reading an early version of the manuscript, pointing out some errors and helping with the proof of \thmref{cluster-cri}.
This work was partially supported by a grant from the National Science Foundation (DMS-0603890) and a grant from the Office of Naval Research (N00014-09-1-0258).

{\small
\bibliographystyle{abbrv}
\bibliography{tree-percolation}
}

\end{document}